\documentclass{article}

\usepackage{graphicx}
\usepackage{amsfonts}
\usepackage{amsmath}
\usepackage{amssymb}

\newtheorem{theorem}{Theorem}

\newtheorem{corollary}{Corollary}

\newtheorem{definition}{Definition}
\newtheorem{example}{Example}

\newtheorem{lemma}{Lemma}

\newtheorem{proposition}{Proposition}
\newtheorem{remark}{Remark}

\newenvironment{proof}[1][Proof]{\textbf{#1.} }{\ \rule{0.5em}{0.5em}}

\begin{document}

\title{On the invariants of some solvable rigid Lie algebras}

\author{Rutwig Campoamor-Stursberg\\Depto. Geometr\'{\i}a y Topolog\'{\i}a\\Fac. CC. Matem\'aticas U.C.M.\\
E-28040 Madrid ( Spain)\newline e-mail: rutwig@nfssrv.mat.ucm.es}

\date{}

\maketitle

\begin{abstract}
We determine fundamental systems of invariants for complex solvable rigid Lie algebras having nonsplit nilradicals of characteristic sequence $\left(3,1,..,1\right)$, these algebras being the natural followers of solvable algebras having Heisenberg nilradicals. A special case of this allows us to obtain a criterion to determine the number of functionally independent invariants of rank one subalgebras of (real or complex) solvable Lie algebras. Finally, we give examples of the inverse procedure, obtaining fundamental systems of an algebra starting from rank one subalgebras.
\end{abstract}

\bigskip

$2000$ \textit{Mathematics Subject Classification: 17B10, 17B30}


\newpage

\section{Introduction}

An important problem arising in the theory of representations of Lie algebras and various physical applications is the determination of invariant functions for the coadjoint representation. Reductions like the Levi decomposition theorem simplify the question to the classes of semi-simple and solvable Lie algebras. The invariants in the semi-simple case were determined in $1950$ by Racah \cite{Ra}. Here the invariants arise as elements of the enveloping algebra, and their number coincides with the rank of the algebra. The study of the non semi-simple case is physically motivated by considerations like internal symmetries of particles or the invariant operators of symmetry groups of physical systems \cite{Ge, Mi}. Important groups as the Galilei and Poincar\'e have been deeply studied \cite{Le, Ma, Mi, Ro}, as well as the subgroups of the latter. These invariants allow to characterize certain systems by giving their energy spectra, angular momenta, etc. For solvable Lie algebras, much less is known. Their invariants have been determined up to dimension six \cite{Nd, Pa}. The invariants found here need not to be polynomials any more, which suggests to call them generalized invariants. The nonexistence of classifications of solvable algebras in dimension $n\geq 7$ force to restrict to concrete classes. In this frame, the generalized Casimir invariants of algebras having an abelian or an Heisenberg Lie algebra as nilradical have recently been computed \cite{NW, Wi}.\newline Starting from this situation, we determine fundamental systems of invariants for all solvable rigid Lie algebras whose nilradical is nonsplit of characteristic sequence $\left(3,1,..,1\right)$. These algebras are the natural followers of those having Heisenberg nilradicals. It is shown that the number of functionally independent invariants is a linear function of their rank. As applications, we deduce a general result on the cardinal of fundamental systems of invariants of rank one subalgebras of a solvable Lie algebra. The converse of this procedure is also of interest, since it allows to calculate fundamental systems of an algebra analyzing certain rank one subalgebras.

\bigskip

We convene that nonwritten brackets are either zero or obtained by antisymmetry. We also use Einstein's convention for sums. Unless otherwise stated, any Lie algebra is nonsplit and complex. 

\section{Generalities}

Let $G$ be a connected Lie group and $\frak{g}$ its Lie algebra. As known, the coadjoint representation is given by 
\[
ad^{*}:G\rightarrow GL(\frak{g}^{*}):\quad (ad_{g}^{*}x)(y)=x(ad_{g^{-1}})y,\quad x,y\in\frak{g}^{*},g\in G
\]
and a function $H\in C^{\infty}\frak{g}^{*}$ is called invariant if $H(x)=H(ad^{*}_{g}x)$. The purpose is to determine a maximal set of functionally independent invariants, which is called a fundamental set of invariants for the coadjoint representation. Probably the most popular procedure to calculate these invariants is the interpretation of the problem in terms of partial differential equations \cite{Le, Ro}. Let $\mathcal{B}=\left\{X_{1},..,X_{n}\right\}$ be a basis of $\frak{g}$ and consider differentiable functions $F\left(x_{1},..,x_{n}\right)$, where the $x_{i}$ are commuting variables on $\frak{g}^{*}$ corresponding to the dual basis $\mathcal{B}$. If $\left\{C_{ij}^{k}\right\}$ are the structure constants of $\frak{g}$ over the preceding basis, then we can find differential operators 
\begin{eqnarray}
\widehat{X_{i}}=-C_{ij}^{k}x_{k}\frac{\partial}{\partial{x_{j}}},\quad 1\leq i\leq n
\end{eqnarray}
realizing the coadjoint representation. The invariants are then obtained from the linear first order partial differential equations
\begin{eqnarray}
\left\{
\widehat{X_{i}}.F=0,\quad 1\leq i\leq n\right.
\end{eqnarray}
The invariants of the coadjoint representation will also be called generalized Casimir invariants, since polynomial invariants give rise to elements in the center of the universal enveloping algebra $\frak{A}$, i.e., to the well known Casimir operators. As follows from the work of Beltrametti and Blasi \cite{Be}, the number of functionally independent invariants of the coadjoint representation $ad^{*}$ of a Lie algebra $\frak{g}$ is given by $\mathcal{N}\left(\frak{g}\right)=dim\left(\frak{g}\right)-r\left(\frak{g}\right)$, where 
$r\left(\frak{g}\right)$ is the maximum rank of the commutator table considered as a matrix.\newline As told before, invariants of solvable Lie algebras have been determined only in low dimensions, due to the impossibility of obtaining classifications in dimensions $n\geq 7$. Two types of solvable algebras have been analyzed in detail in arbitrary dimension: those having abelian nilradical and those having it isomorphic to the Heisenberg Lie algebra $\frak{h}_{k}$. From the latter we easily see that the direct sum of $\frak{h}_{k}$ and a maximal torus of derivations is rigid. Since for rigid Lie algebras the torus determines completely the isomorphism class of the nilradical and these algebras decompose in a particular manner, in addition to the existing algorithms that allow (theoretically) a classification of solvable rigid Lie algebras \cite{AG2}, it is reasonable to focus on this class and their invariants. We therefore recall the most elementary facts about rigidity: Let $\frak{L}^{n}$ be the algebraic variety of Lie algebras laws of
dimension $n$ over $\mathbb{C}$. If $\mu \in \frak{L}^{n}$ and $\left\{
e_{i}\right\} _{1\leq i\leq n}$ is a basis of $\mathbb{C}^{n},$ then the
structure constants of the law $\mu $ satisfy the following polynomial
equations : 
\begin{eqnarray}
C_{ij}^{k}+C_{ji}^{k}=0,\;\;1\leq i<j\leq n,\;1\leq k\leq n \\ 
\sum_{l=1}^{n}C_{ij}^{l}C_{lk}^{s}+C_{jk}^{l}C_{li}^{s}+C_{ki}^{l}C_{jl}^{s}=0,\quad 1\leq i<j<k\leq n,\;1\leq s\leq n
\end{eqnarray}

So the $C_{ij}^{k}$ satisfying the system above parametrize the variety $%
\frak{L}^{n},$ embedded in $\mathbb{C}^{n^{3}}$. The subset $\frak{N}^{n}$ of nilpotent Lie algebras is closed in this variety. From now on, we identify a Lie algebra with its law. The linear group $GL\left( n,\mathbb{C}\right) $
acts on $\frak{L}^{n}$ by: 
\[
GL\left( n,\mathbb{C}\right) \times \frak{L}^{n}\rightarrow \frak{L}^{n},\quad \left( f,\mu \right) \mapsto f^{-1}\mu \left( f,f\right) 
\]
where 
\[
\mu \left( f,f\right) \left( X,Y\right) =\mu \left( f\left( X\right)
,f\left( Y\right) \right) \quad \forall X,Y\in \mathbb{C}^{n} 
\]
The resulting orbit of the law $\mu $ 
\begin{eqnarray}
O\left( \mu \right) =\left\{ h^{-1}\mu \left( h,h\right) \;|\;h\in GL\left(
n,\Bbb{C}\right) \right\} 
\end{eqnarray}
consists of all laws isomorphic to $\mu $.

\begin{definition}
A law $\mu \in \frak{L}^{n}$ is called rigid if its orbit 
$O\left( \mu \right) $ is open in $\frak{L}^{n}$.
\end{definition}

Let $\mu \in \frak{L}^{n}$ be a solvable, non nilpotent law. The
corresponding Lie algebra $\frak{g}$ is called decomposable if $\frak{g}$
admits a decomposition 
\[
\frak{g}=\frak{n}\oplus \frak{t} 
\]
where $\frak{n}$ is the nilradical of $\frak{g}$ and $\frak{t}$ is an
abelian subalgebra consisting of $ad-$semisimple elements. $\frak{t}$ is
called exterior torus of derivations. The following theorem is an immediate
consequence of Carles algebraicity theorem \cite{Ca} :

\begin{theorem}
\textit{A rigid solvable Lie algebra }$\frak{g}$ is decomposable.
\end{theorem}

In order to describe this decomposition properly, we have to find adequate elements in the torus $\frak{t}$ that describe the action on the nilradical. Let $X\in \frak{t}$ be a nonzero vector. We call $X$ regular if the
dimension of the kernel of the adjoint operator is minimal among the
elements of $\frak{t}$, i.e 
\[
\dim \;Ker\,ad\left( X\right) =\min \;\left\{ \dim \;Ker\,ad\left( Y\right)
\;|\;Y\in \frak{t}\right\} 
\]
If $q=\dim \;Ker\,ad(X),$ then we can find a basis formed by eigenvectors of
the adjoint operator $\left( X_{1},...,X_{n}=X\right) $ such that $\left(
X_{1},..X_{p+q}\right) $ is a basis of the nilradical $\frak{n}$, $\left(
X_{p+1+q},..,X_{n}\right) $ is a basis of $\frak{t}$ and $\left(
X_{p+1},..,X_{n}\right) $ is a basis of $Ker\;ad\left( X\right) .$\newline
The linear system of roots associated to $\left( X_{1},..,X_{n}\right) $ is
the linear system to $n-1$ variables $x_{i}$ whose equations are 
\begin{eqnarray}
x_{i}+x_{j}=x_{k} 
\end{eqnarray}
if the component of the vector $\mu \left( X_{i},X_{j}\right) $ on $X_{k}$
is nonzero. We will note this system by $S\left( X\right) $ or simply $S.$
The following result establishes a necessary condition for a Lie algebra to
be rigid \cite{AG1}:

\begin{theorem}
( of the rank ). \textit{If the solvable, decomposable Lie algebra }$\frak{g}%
=\frak{n}\oplus \frak{t}$ \textit{is rigid, then for all regular vectors }$X$
\textit{we have} 
\[
rank\;S\left( X\right) =\dim \,\frak{n}-1
\]
\textit{In particular the rank does not depend on the choice of the basis or
the regular vector.}
\end{theorem}

\begin{corollary}
\textit{If }$\frak{g}=\frak{n}\oplus \frak{t}$ \textit{is rigid, then }$%
\frak{t}$ \textit{is a maximal torus.}
\end{corollary}

A proof of the preceding result can be found in \cite{AG1}. It allows to search
rigid laws starting with a sequence of weights for a fixed regular vector $X$. For such an algebra, we note the rank by $rg\left( \frak{g}\right) =\dim \frak{t}$.
\newline We finally recall an invariant of nilpotent Lie algebras that will be helpful: Let $\frak{g}$ be a nilpotent Lie algebra and $X\in \frak{g}-C^{1}\left( 
\frak{g}\right) $ a nonzero vector. We denote by $c\left( X\right) $ the ordered sequence of dimensions of the
Jordan blocks for the adjoint operator. Considering the lexicographical order in
the set of these sequences 
\[
\left( c_{1},..,c_{t}\right) \geq \left( c_{1}^{\prime },..,c_{s}^{\prime
}\right) \leftrightarrow \exists i \textrm{ such that } c_{j}=c_{j}^{\prime }%
\textrm{ for } i>j \textrm{ and }c_{i}>c_{i}^{\prime } 
\]
we define the characteristic sequence of $\frak{g}$ as 
\[
c\left( \frak{g}\right) =\max_{X\in \frak{g}-C^{1}\left( \frak{g}\right)
}\left\{ c\left( X\right) \right\} 
\]
A vector $X\neq 0,X\in \frak{g}-C^{1}\frak{g}$ satisfying $c\left( \frak{g}%
\right) =c\left( X\right) $ is called characteristic vector.

\begin{remark}
It is obvious that abelian Lie algebras have characteristic sequence $\left(1,1,..,1\right)$, while for the Heisenberg Lie algebra $\frak{h}_{k}$ we obtain $c(\frak{h}_{k})=\left(2,1,..,1\right)$. Thus it is reasonable to study the next step, the sequence $\left(3,1,..,1\right)$ and to determine the rigid Lie algebras having a nilradical with this sequence.
\end{remark}

\begin{theorem}
\textrm{\cite{AC1}} Any solvable rigid law whose nilradical $\frak{n}$ is nonsplit and of characteristic sequence $\left(3,1,..,1\right)$ is isomorphic to one of the following laws:

\begin{enumerate}
\item  The Lie algebra $\frak{d}_{2m}=\frak{g}_{2m}\oplus \frak{t}\quad (m\geq 3)$
with basis $(X_{0},X_{1},X_{2},X_{3},Y_{1},..,Y_{2m-4},$\newline
$V_{1},..V_{m})$ and law 

\begin{eqnarray*}
\lbrack X_{0},X_{i}] &=&X_{i+1},\;i=1,2 \\
\lbrack Y_{2i-1},Y_{2i}] &=&X_{3},\;1\leq i\leq m-2 \\
\lbrack V_{1},X_{i}] &=&\left( i+1\right) X_{i},\;i=0,1,2,3 \\
\lbrack V_{1},Y_{2i-1}] &=&\left( i+4\right) Y_{2i-1},\;1\leq i\leq m-2 \\
\lbrack V_{1},Y_{2i}] &=&-i.Y_{2i},\;1\leq i\leq m-2 \\
\lbrack V_{2},X_{i}] &=&X_{i},\;i=1,2,3 \\
\lbrack V_{2},Y_{2i}] &=&Y_{2i},\;1\leq i\leq m-2 \\
\lbrack V_{i+2},Y_{2i-1}] &=&Y_{2i-1},\;1\leq i\leq m-2 \\
\lbrack V_{i+2},Y_{2i}] &=&-Y_{2i},\;1\leq i\leq m-2
\end{eqnarray*}

\item  The Lie algebra $\frak{d}_{2m+1}=\frak{g}_{2m+1}\oplus 
\frak{t}\quad (m\geq 2)$ with basis $(X_{0},..,X_{3},Y_{1},..,Y_{2m-3},$\newline
$V_{1},..V_{m}$ and law
 
\begin{eqnarray*}
\lbrack X_{0},X_{i}] &=&X_{i+1},\;i=1,2 \\
\lbrack X_{1},Y_{2m-3}] &=&X_{3} \\
\lbrack Y_{2i-1},Y_{2i}] &=&X_{3},\;1\leq i\leq m-2 \\
\lbrack V_{1},X_{i}] &=&\left( i+1\right) X_{i},\;i=0,1,2,3 \\
\lbrack V_{1},Y_{2i-1}] &=&\left( i+4\right) Y_{2i-1},\;1\leq i\leq m-2 \\
\lbrack V_{1},Y_{2i}] &=&-iY_{2i},\;1\leq i\leq m-2 \\
\lbrack V_{1},Y_{2m-3}] &=&2Y_{2m-3}\\
\lbrack V_{2},X_{i}] &=&X_{i},\;i=1,2,3 \\
\lbrack V_{2},Y_{2i}] &=&Y_{2i},\;1\leq i\leq m-2 \\
\lbrack V_{i+2},Y_{2i-1}] &=&Y_{2i-1},\;1\leq i\leq m-2 \\
\lbrack V_{i+2},Y_{2i}] &=&-Y_{2i},\;1\leq i\leq m-3
\end{eqnarray*}

\item  The 7 dimensional Lie algebra $\frak{d}_{5}^{\prime }$
 with basis $\left( X_{0},..,X_{3},Y_{1},V_{1},V_{2}\right) $%
 and law 
 
\begin{eqnarray*}
\lbrack X_{0},X_{i}] &=&X_{i+1},\;i=1,2 \\
\lbrack X_{1},X_{2}] &=&Y_{1} \\
\lbrack V_{1},X_{i}] &=&\left( i+1\right) X_{i},\;i=0,1,2,3 \\
\lbrack V_{1},Y_{1}] &=&5Y_{1} \\
\lbrack V_{2},X_{i}] &=&X_{i},\;i=1,2,3 \\
\lbrack V_{2},Y_{1}] &=&2Y_{1}
\end{eqnarray*}
\end{enumerate}
\end{theorem}

\begin{remark} Indeed more is true: any nonsplit nilpotent Lie algebra $\frak{n}$ of characteristic sequence $\left(3,1,..,1\right)$ is isomorphic to the nilradical of one of the former solvable algebras \cite{AC1}. 
\end{remark}

In order to determine the invariants of these algebras, we have to find the cardinal of their fundamental systems of invariants: 

\begin{lemma}
Following identities hold:
\begin{enumerate}

\item For $m\geq 3$ we have $\mathcal{N}(\frak{d_{2m}})=m-2$.

\item For $m\geq 2$ we have $\mathcal{N}(\frak{d_{2m+1}})=m-1$.

\item $\mathcal{N}(\frak{d}_{5}^{\prime})=1$.
\end{enumerate}
\end{lemma}

\begin{proof}
We prove the result for $\frak{d}_{2m+1}$, the reasoning for $\frak{d}_{2m}$ being similar. For $\frak{d}_{5}^{\prime}$ the result is trivial. Clearly the nilradical of $\frak{d}_{2m+1}$ is isomorphic to the algebra spanned by $\left\{X_{0},X_{1},Y_{i}\right\}_{1\leq i\leq m-3}$ and the torus $\frak{t}$ is generated by the $\left\{V_{i}\right\}$. The root system $S$ associated to $\frak{d}_{2m+1}$ is
\begin{eqnarray}
\lambda _{0}+\lambda _{1}=\lambda _{2} \\ 
2\lambda _{0}+\lambda _{1}=\lambda _{3} \\
\mu_{2i-1}+\mu_{2i}=\lambda_{3} \\
\mu _{2m-3}=2\lambda _{0} 
\end{eqnarray}

The matrix $A=\left(C_{ij}^{k}\right)$ of structure constants for $\frak{d}_{2m+1}$ over this basis can be extracted from the following table

\bigskip  

$
\begin{tabular}{|ccccc|cc|cc|c|c|}
\hline\hline
$\mathbf{x}_{0}$ & \multicolumn{1}{|c}{$\mathbf{x}_{1}$} & 
\multicolumn{1}{|c}{$\mathbf{x}_{2}$} & \multicolumn{1}{|c}{$\mathbf{x}_{3}$}
& \multicolumn{1}{|c|}{$\mathbf{y}_{2m-3}$} & \multicolumn{1}{||c}{$\mathbf{y%
}_{2i-1}$} & \multicolumn{1}{|c|}{$\mathbf{y}_{2i}$} & \multicolumn{1}{||c}{$%
\mathbf{v}_{1}$} & \multicolumn{1}{|c|}{$\mathbf{v}_{2}$} & $\mathbf{v}_{i+2}
$ & $[]$ \\ \hline\hline
$0$ & $x_{2}$ & $x_{3}$ & $0$ & $0$ & $0$ & $0$ & $-x_{0}$ & $0$ & $0$ & 
\multicolumn{1}{||c||}{$\mathbf{x}_{0}$} \\ \cline{11-11}
$-x_{2}$ & $0$ & $0$ & $0$ & $x_{3}$ & $0$ & $0$ & $0$ & $-x_{1}$ & $0$ & 
\multicolumn{1}{||c||}{$\mathbf{x}_{1}$} \\ \cline{11-11}
$-x_{3}$ & $0$ & $0$ & $0$ & $0$ & $0$ & $0$ & $-x_{2}$ & $-x_{2}$ & $0$ & 
\multicolumn{1}{||c||}{$\mathbf{x}_{2}$} \\ \cline{11-11}
$0$ & $0$ & $0$ & $0$ & $0$ & $0$ & $0$ & $-2x_{3}$ & $-x_{3}$ & $0$ & 
\multicolumn{1}{||c||}{$\mathbf{x}_{3}$} \\ \cline{11-11}
$0$ & $-x_{3}$ & $0$ & $0$ & $0$ & $0$ & $0$ & $-2y_{2m-3}$ & $0$ & $0$ & 
\multicolumn{1}{||c||}{$\mathbf{y}_{2m-3}$} \\ \hline
$0$ & $0$ & $0$ & $0$ & $0$ & $0$ & $-x_{3}$ & $0$ & $0$ & $-y_{2i-1}$ & 
\multicolumn{1}{||c||}{$\mathbf{y}_{2i-1}$} \\ \cline{11-11}
$0$ & $0$ & $0$ & $0$ & $0$ & $x_{3}$ & $0$ & $-2y_{2i}$ & $-y_{2i}$ & $%
y_{2i}$ & \multicolumn{1}{||c||}{$\mathbf{y}_{2i}$} \\ \hline
$x_{0}$ & $0$ & $x_{2}$ & $2x_{3}$ & $2y_{2m-3}$ & $0$ & $2y_{2i}$ & $0$ & $0
$ & $0$ & \multicolumn{1}{||c||}{$\mathbf{v}_{1}$} \\ \cline{11-11}
$0$ & $x_{1}$ & $x_{2}$ & $x_{3}$ & $0$ & $0$ & $y_{2i}$ & $0$ & $0$ & $0$ & 
\multicolumn{1}{||c||}{$\mathbf{v}_{2}$} \\ \hline
$0$ & $0$ & $0$ & $0$ & $0$ & $y_{2i-1}$ & $-y_{2i}$ & $0$ & $0$ & $0$ & 
\multicolumn{1}{||c||}{$\mathbf{v}_{i+2}$} \\ \hline
\end{tabular}
$

\bigskip
 
where the columns between double lines are representative for the vectors  $\left\{y_{2i-1},y_{2i}\right\},\quad 1\leq i\leq m-2$. For any $1\leq i,j\leq 3m+1$ let $A_{\left( i,j\right) }$ be the matrix
obtained by deletion of the $i^{th}$ row and $j^{th}$ column of $A$. By
recurrence, let 
\[
A_{\left( i_{1},j_{1}\right) ,..,\left( i_{r},j_{r}\right) }=\left(
A_{\left( i_{1},j_{1}\right) ,..,\left( i_{r-1},j_{r-1}\right) }\right)
_{\left( i_{r},j_{r}\right) }
\]
For $m\geq 2$ define 
\[
A_{m}:=A_{\left( 3m+1,3m+1\right) ,..,\left( 2m+3,2m+3\right) }
\]
It can easily be proven by induction over $m$ that $\det \left(
A_{m+1}\right) =x_{3}^{2}\det \left( A_{m}\right) $, and since $\det \left(
A_{2}\right) =4x_{3}^{6}$ we get $\det \left( A_{m}\right) =4x_{3}^{2m+2}$.
From this we deduce $rank\left( A\right) =2m+2$. Since the vector $V_{1}$ is regular, the root system $(S)$ describes the torus action and the rank of the algebra. As a consequence we have that 
\[
2m+2=\sup \left\{ rank\,A\left( \mu \right) \;|\;\mu \in \mathcal{O}\left(
\frak{g}_{2m+1}\right) \right\} 
\]
and by the Beltrametti-Blasi formula we obtain $\mathcal{N}(\frak{d}_{2m+1})=m-1$.
\end{proof}

\begin{theorem}
Let $\frak{r}$ be a solvable rigid Lie algebra whose nilradical is nonsplit
 of characteristic sequence $\left( 3,1,..,1\right) $. Then a set of
fuctionally independent invariants is given by 

\begin{enumerate}
\item  $\left\{ \frac{y_{2k-1}y_{2k}+x_{3}v_{k+2}}{x_{3}}\right\} _{1\leq
k\leq m-2}$ if $\frak{r}\simeq \frak{d}_{2m}$

\item  $\left\{ \frac{%
x_{0}x_{2}x_{3}+x_{2}^{2}y_{2m-3}-2x_{3}^{2}v_{2}+v_{1}x_{3}^{2}-2x_{1}x_{3}y_{2m-3}%
}{-2x_{3}^{2}},\frac{x_{3}v_{k+2}+y_{2k-1}y_{2k}}{x_{3}}\right\}_{1\leq k\leq m-2}$ if $\frak{r}\simeq \frak{d}_{2m+1}$

\item  $\left\{ \frac{\left( 2x_{0}y_{1}+x_{2}^{2}-2x_{1}x_{3}\right) ^{3}}{%
\left( x_{3}y_{1}\right) ^{2}}\right\} $ if $\frak{r}\simeq \frak{d}_{5}^{\prime}$
\end{enumerate}
\end{theorem}

\begin{proof}

\begin{enumerate}

\item Let $\frak{r}=\frak{d}_{2m+1}$. The system $(2.2)$ is in this case: 

\begin{eqnarray}
\widehat{X_{0}}=-x_{2}\partial _{x_{1}}-x_{3}\partial
_{x_{2}}+x_{0}\partial _{v_{1}}\\
\widehat{X_{1}}=x_{2}\partial _{x_{0}}+x_{1}\partial _{v_{2}}\\
\widehat{X_{2}}=x_{3}\partial _{x_{0}}+x_{2}\partial _{v_{1}}+x_{2}\partial
_{v_{2}}\\
\widehat{X_{3}}=2x_{3}\partial _{v_{1}}+x_{3}\partial _{v_{2}}\\
\widehat{Y_{2i-1}}=-x_{3}\partial _{y_{2i}}+y_{2i-1}\partial
_{v_{i+1}},\quad 1\leq i\leq m-2\\
\widehat{Y_{2i}}=x_{3}\partial _{y_{2i-1}}-y_{2i}\partial
_{v_{i+1}}+2y_{2i}\partial _{v_{1}}+y_{2i-1}\partial _{v_{2}},\quad 1\leq i\leq m-2\\
\widehat{V_{1}}=x_{0}\partial _{x_{0}}+x_{2}\partial
_{x_{2}}+2x_{3}\partial _{x_{3}}+2\sum_{i=1}^{m-2}y_{2i}\partial _{y_{2i}}\\
\widehat{V_{2}}=x_{1}\partial _{x_{1}}+x_{2}\partial _{x_{2}}+x_{3}\partial
_{x_{3}}+\sum_{i=1}^{m-2}y_{2i}\partial _{y_{2i}}\\
\widehat{V_{i+1}}=y_{2i-1}\partial _{y_{2i-1}}-y_{2i}\partial _{y_{2i}},\quad 1\leq i\leq m-2
\end{eqnarray}

Elementary manipulations show that an invariant $F$ satisfies $\partial_{x_{i}}F=0$ for $i=0,1,2$ and $\partial_{v_{i}}F=0$ for $i=1,2$. This reduces the system to the following 
\begin{eqnarray}
\left( -x_{3}\partial _{y_{2i}}+y_{2i-1}\partial _{v_{i+1}}\right) F=0,1\leq
i\leq m-2 \\ 
\left( x_{3}\partial _{y_{2i-1}}-y_{2i}\partial _{v_{i+1}}+2y_{2i}\partial
_{v_{1}}+y_{2i-1}\partial _{v_{2}}\right) F=0,1\leq i\leq m-2 \\ 
\left( x_{3}\partial _{x_{3}}+\sum_{i=1}^{m-2}y_{2i}\partial
_{y_{2i}}\right) F=0 \\ 
\left( y_{2i-1}\partial _{y_{2i-1}}-y_{2i}\partial _{y_{2i}}\right) F=0
\end{eqnarray}

This system, which can be integrated by elementary methods \cite{Ol}, gives a  fundamental set of invariants formed by the functions
\[
\left\{f_{i}=\frac{y_{2i-1}y_{2i}+x_{3}v_{i+2}}{x_{3}}\right\}_{1\leq i\leq m-2}
\]

\item Let $\frak{r}=\frak{d}_{2m+1}$. Consider the following system of PDE's :%
\begin{eqnarray}
\widehat{X_{i}}.F=0,\;0\leq i\leq 3 \\ 
\widehat{Y_{2m-3}}.F=0
\end{eqnarray}

It corresponds to the system of PDEs that we obtain for $m=2$, and by the preceding lemma it has one solution, which can be determined by elementary
methods. We obtain 
\[
f_{1}=\frac{%
2x_{3}^{2}v_{2}-x_{3}^{2}v_{1}+2x_{1}x_{3}y_{2m-3}-x_{0}x_{2}x_{3}-x_{2}^{2}y_{2m-3}%
}{x_{3}^{2}}
\]
Thus, for $m=2$, a fundamental set of invariants is given by $\left\{ f_{1}\right\} $. For the
general case $m\geq 3$ observe that, since $\widehat{X_{3}}.F=x_{3}\left(
2\partial _{v_{1}}+\partial _{v_{2}}\right) =0$, the equation $\widehat{%
Y_{2i}}.F=0$ can be reduced to 
\begin{eqnarray}
\widehat{Y_{2i}^{\prime }}.F=\left( x_{3}\partial
_{y_{2i-1}}-y_{2i}\partial _{v_{i+2}}\right) =0
\end{eqnarray}
As a consequence, the function $f_{1}$ is a solution of $\left( 2.25-26\right) $
for any $m\geq 3$, and the remaining fundamental solutions are obtained from
the system 
\begin{eqnarray}
\widehat{Y_{2i-1}}.F=0,\;1\leq i\leq m-2 \\ 
\widehat{Y_{2i}^{\prime }}.F=0,\;1\leq i\leq m-2
\end{eqnarray}

For any fixed $i\in \left\{ 1,..,m-2\right\} $ a solution of the
corresponding equations is 
\[
f_{i+1}=\frac{x_{3}v_{i+2}+y_{2i-1}y_{2i}}{x_{3}}
\]
which shows that a fundamental system of invariants for $\frak{d}_{2m+1}$ is given by $\left\{
f_{1},f_{i+1}\right\} _{1\leq i\leq m-2}$. 

\item It is easily verified that any invariant satisfies $\partial_{v_{1}}F=\partial_{v_{2}}F=0$. Transformations of the corresponding system show that an invariant function $F$ must be a
solution of the following equation :%
\begin{eqnarray}
\left( -x_{2}^{2}+y_{1}x_{0}+2x_{1}x_{3}\right) \partial
_{x_{1}}F+3y_{1}x_{3}\partial _{y_{1}}F=0
\end{eqnarray}
This equation admits, among others, the rational solution 
\begin{eqnarray}
F_{1}=\frac{\left( 2x_{0}y_{1}+x_{2}^{2}-2x_{1}x_{3}\right) ^{3}}{%
x_{3}^{2}y_{1}^{2}}
\end{eqnarray}
and it can easily verified that this function satisfies the system. Thus a
set of invariants is given by $\left\{ F_{1}\right\} $.
\end{enumerate}
\end{proof}

\section{Applications to rank one subalgebras}

The seven dimensional rigid Lie algebra $\frak{d}_{5}^{\prime }$ is a nice
example to illustrate how, under certain assumptions, we can deduce the
invariants of some subalgebras. Let us consider the rank one solvable
(on-nilpotent) subalgebras $\frak{d}_{5}^{a,b}$ of $\frak{d}_{5}^{\prime }$
having the following brackets :\newline
\[
\begin{tabular}{lll}
$\left[ X_{0},X_{1}\right] =X_{2},$ & $\left[ X_{0},X_{2}\right] =X_{3},$ & $%
\left[ X_{1},X_{2}\right] =Y_{1},$ \\ 
$\left[ V,X_{0}\right] =aX_{0},$ & $\left[ V,X_{1}\right] =bX_{1},$ & $%
\left[ V,X_{2}\right] =\left( a+b\right) X_{2},$ \\ 
$\left[ V,X_{3}\right] =\left( 2a+b\right) X_{3},$ & $\left[ V,Y_{1}\right]
=\left( a+2b\right) Y_{1};$ & $a,b\in\mathbb{C}$%
\end{tabular}
\]
Clearly we have $V=aV_{1}+bV_{2}$. The system of PDEs giving the invariants
of $\frak{d}_{5}^{a,b}$ is
\begin{eqnarray*}
\widehat{X_{i}}F=0,\;0\leq i\leq 3 \\ 
\widehat{Y_{1}}F=0\\ 
\widehat{V}F=0
\end{eqnarray*}
where 
\begin{eqnarray*}
\widehat{V}=\left(ax_{0}\partial _{x_{0}}+bx_{1}\partial _{x_{1}}+\left(
a+b\right)x_{2}\partial _{x_{2}}+\left(2a+b\right)x_{3}\partial
_{x_{3}}+\left( a+2b\right)y_{1}\partial _{y_{1}}\right)
\end{eqnarray*}

Now any invariant $F$ satisfies $\partial _{v}F=0$, we obtain that any
invariant of $\frak{d}_{5}^{\prime }$ is also an invariant of $\frak{d}%
_{5}^{a,b}$, since the five first equations coincide and the sixth is a
linear combination of the toral equations of $\frak{d}_{5}^{\prime }$. 

\begin{lemma}
For any $a,b\in\mathbb{C}$ not simultaneously zero we have $\mathcal{N}\left( \frak{d}%
_{5}^{\prime }\right) =2$.
\end{lemma}

\begin{proof}
Since $F_{1}$ satisfies the system, $\mathcal{N}\geq 1$. Now it can easily
be seen that the matrix obtained from the commutation table has rank four,
from which the assertion follows.
\end{proof}

We make a distinction among the subalgebras $\frak{d}_{5}^{\prime }$: those
having trivial center and those not. It is a routinary verification that $%
\frak{d}_{5}^{a,b}$ has nontrivial center if and only if $a+2b=0$ or $2a+b=0$%
. Since solvable Lie algebras having nontrivial center admit at least a
polynomial invariant, we obtain the following 

\begin{lemma}
Let $\frak{d}_{5}^{a,b}$ have nontrivial center. Then a fundamental set of
invariants is given by 

\begin{enumerate}
\item  $\left\{ F_{1},x_{3}\right\} $ if $2a+b=0,$

\item  $\left\{ F_{1},y_{1}\right\} $ if $a+2b=0$.
\end{enumerate}
\end{lemma}

For the case of a non-trivial center, it is easy to verify that any invariant of $\frak{d}_{5}^{a,b}$ is a function of $F_{1},x_{3}$ and $%
y_{1}$. Specifically, we have 

\begin{lemma}
If $Z\left(\frak{d}_{5}^{a,b}\right)\neq 0$, a fundamental set of invariants for $\frak{d}_{5}^{a,b}$ is given by $\left\{F_{1}, \frac{y_{1}^{2a+b}}{x_{3}^{a+2b}}\right\}$.
\end{lemma}

The interest of this example is that it provides us with a general result: Let $\frak{r}=\frak{n}\oplus\frak{t}$ be a solvable non-nilpotent Lie algebra over $\mathbb{K}=\lbrace \mathbb{R},\mathbb{C}\rbrace$ and $\frak{t}$ a torus of derivations. Let $\left\{X_{1},..,X_{n},V_{1},..,V_{p}\right\}$ be a basis such that the nilradical $\frak{n}$ is generated by the $X_{i}$ and the torus by the $V_{i}$. Suppose moreover that $\frak{r}$ satisfies the following property:

\begin{eqnarray}
\partial_{v_{i}}F=0,\quad 1\leq i\leq p \textrm{ for any invariant } F
\end{eqnarray}

Let $\frak{t}^{\prime}$ be a one dimensional subtorus of $\frak{t}$ and $\frak{r}^{\prime}=\frak{n}\oplus\frak{t}^{\prime}$ such that $\frak{t}^{\prime}.Z\left(\frak{n}\right)\neq 0$, i.e., the action of the subtorus $\frak{t}^{\prime}$ on the center of the nilradical is nonzero.

\begin{proposition}
In the preceding conditions, any element $F$ of a fundamental set of invariants for $\frak{r}$ is an invariant of $\frak{r}^{\prime}$. Moreover, $\mathcal{N}(\frak{r}^{\prime})\geq \mathcal{N}(\frak{r})$.
\end{proposition}

\begin{proof}
Let $V$ generate $\frak{t}^{\prime}$. We can find coefficients $\alpha_{1},..,\alpha_{p}\in\mathbb{K}$ such that $V=\sum_{i=1}^{p}\alpha_{i}V_{i}$. Since $\frak{n}$ is nilpotent, its center is nonzero, and by $\frak{t}^{\prime}.Z\left(\frak{n}\right)\neq 0$, we deduce that $\partial_{v}F=0$ for any invariant $F$ of $\frak{r}^{\prime}$. The system of PDEs giving the invariants is 

\begin{eqnarray}
\widehat{X_{i}}.F=-\lbrack X_{i},X_{j}\rbrack \partial_{x_{j}}.F=0, & 1\leq i\leq n \\
\widehat{V}.F=-\lbrack V,X_{i}\rbrack \partial_{x_{i}}.F=0 & 
\end{eqnarray}

Now $\widehat{V}=\sum_{i=1}^{p}\alpha_{i}\widehat{V_{i}}$. This shows that any invariant $G$ of $\frak{r}$ also satisfies the system $(3.6),(3.7)$, since the equations concerning the operators $\widehat{X_{i}}$ coincide, and the last equation is a linear combination of the equations $\lbrace \widehat{V_{i}}.F=0 \rbrace$ of $\frak{r}$. Thus any fundamental set of invariants for $\frak{r}$ is a set of functionally independent invariants for $\frak{r}^{\prime}$. As a consequence $\mathcal{N}(\frak{r}^{\prime})\geq \mathcal{N}(\frak{r})$.
\end{proof}

\begin{remark}
If $p$ is even, by the Beltrametti-Blasi formula we deduce  $\mathcal{N}(\frak{r}^{\prime})\geq \mathcal{N}(\frak{r})+1$
\end{remark}

The main interest of this result is its possible application in the other direction [whenever it makes sense], i.e., starting from the invariants of a rank one subalgebra to find a fundamental set of the algebra $\frak{r}$. We illustrate it by an example:

\begin{example}
Consider the solvable Lie algebra $\frak{g}$ given by the brackets 
\begin{eqnarray*}
\lbrack V_{1},Y_{i}\rbrack = iY_{i},\quad 1\leq i\leq 7\\
\lbrack V_{2},Y_{i}\rbrack = y_{i},\quad 2\leq i\leq 7\\
\lbrack Y_{1},Y_{i}\rbrack = Y_{i+1},\quad 2\leq i\leq 6\\
\end{eqnarray*}
This algebra is of rank $2$, and a torus is generated by $V_{1}$ an $V_{2}$. From the corresponding system (2.2) it can be seen that for any invariant $F$ of $\frak{g}$ we have $\partial_{v_{i}}F=0, i=1,2$, thus we can apply the preceding proposition. Consider the one dimensional subtorus generated by $V_{2}$. The semidirect product of this torus and the nilradical is an eight dimensional Lie algebra $\frak{g}^{\prime}$ satisfying $\mathcal{N}(\frak{g}^{\prime})=2$. A fundamental set of invariants for this algebra is formed by the following functions:
\begin{eqnarray*}
g_{1}=\frac{6y_{4}y_{6}y_{7}^{2}-6y_{3}y_{7}^{3}+2y_{6}^{4}-8y_{5}y_{6}^{2}y_{7}+5y_{5}^{2}y_{7}^{2}}{3y_{4}y_{7}^{3}+y_{6}^{3}y_{7}-3y_{5}y_{6}y_{7}^{2}}\\
g_{2}=\frac{-5y_{2}y_{7}^{4}+5y_{4}y_{5}y_{7}^{3}+5y_{5}y_{6}^{3}y_{7}-5y_{5}^{2}y_{6}y{7}^{2}-y_{6}^{5}+5y_{3}y_{6}y_{7}^{3}-5y_{4}y_{6}^{2}y_{7}^{2}}{3y_{4}y_{7}^{4}+y_{6}^{3}y_{7}^{2}-3y_{5}y_{6}y_{7}^{3}}\\
\end{eqnarray*}
As the assumption of the proposition are satisfied, we know that there exists a function $F(g_{1},g_{2})$ that forms a fundamental system of invariants for the Lie algebra $\frak{g}$. In this case it is not difficult to see that such a function is given by $F=\frac{g_{2}}{g_{1}^{2}}$. 
\end{example}

\begin{remark}
Of course this inverse procedure is applicable only in quite concrete cases, since the algebra $\frak{r}$ could have no invariants, while a rank one odd-dimensional Lie subalgebra necessarily has solutions.
\end{remark}

\section*{Concluding remarks}

The preceding results show that rigid Lie algebras are an adequate class for studying the generalized Casimir invariants. Since the torus determines the nilradical, the equations describing the action of the torus are of importance for the structure of the invariants, specifically if those invariants do not depend on the variables corresponding to toral elements. Since rigid Lie algebras can theoretically be classified following the algorithm given in \cite{AG2}, the analysis these algebras in relatively low dimension and their subalgebras can provide alternative criteria to attack the general case, as classifications in general are not possible. This in particular applies to the class of $2$-step solvable Lie algebras, for which general constructions exist \cite{C3}. Observe further that, like proposition 1, the results do not depend on the field, which allows to derive results for both real and complex Lie algebras. The reason to analyze the complex rigid laws lies in their decomposition: for real rigid Lie algebras no result like theorem 1 is known, and it is still an open question wheter it can exist.\newline Finally, the same analysis could be applied to any solvable Lie algebra whose nilradical has characteristic sequence $\left(p,1,..,1\right)$ for $p\geq 4$. These algebras, which are completely classified up to $p=5$ \cite{AC2, C2}, consist of a low dimensional subalgebra $\frak{m}$ generated by a characteristic vector, to which $3$-dimensional Heisenberg Lie algebras having its derived subalgebra in the center of $\frak{m}$ are glued. This structure makes this class worthy of being analyzed to obtain alternative criteria for the generalized Casimir invariants of solvable Lie algebras.

\end{document}